\documentclass[amsthm]{elsart}
\usepackage{yjsco}
\usepackage{natbib}

\usepackage{amssymb}
\usepackage{graphicx}%
\newenvironment{undef}[1]%
           {\vspace{3.3mm}
           \noindent{\bf #1}\it}%
           {\vspace{3.3mm}}

\def\Res{{\rm{Res}}}

\def\Sres{{\mbox{Sres}}}

\def\V{{\mathcal V}}

\def\II{{\mathbb I}}

\def\N{{\mathbb N}}

\def\Res{{\rm{Res}}}

\begin{document}

\begin{frontmatter}

\title{Sylvester's Double Sums: the general case}

\thanks{D'Andrea's Research supported by the Research Project MTM2007-67493 of the Ministerio de Educaci\'on y Ciencia, Spain;
Hong's Research supported by NSF CCR-0097976; Krick's research supported by Argentine grants ANPCyT PICT  33671/2005 and UBACYT X-112;  Szanto's
research supported by NSF grants CCR-0306406 and CCR-0347506}
\author{Carlos D'Andrea}
\address{Department d'\`Algebra i Geometria, Facultat de
Matem\`atiques, Universitat de Barcelona, Gran Via de les Corts
Catalanes, 585; 08007 Spain.} \ead{cdandrea@ub.edu}
\ead[url]{http://carlos.dandrea.name}

\author{Hoon Hong}
\address{Department of Mathematics, North Carolina State University,
Raleigh NC 27695, USA.}
\ead{hong@math.ncsu.edu}
\ead[url]{http://www4.ncsu.edu/\~\,hong}

\author{Teresa Krick}
\address{Departamento de Matem\'atica, Facultad de Ciencias Exactas y
Naturales, Universidad de Buenos Aires, Ciudad Universitaria, 1428
Buenos Aires, Argentina and CONICET, Argentina.}
\ead{krick@dm.uba.ar} \ead[url]{http://mate.dm.uba.ar/\~\,krick}

\author{Agnes Szanto}
\address{Department of Mathematics, North Carolina State
University, Raleigh NC 27695, USA.}
\ead{aszanto@ncsu.edu}
\ead[url]{http://www4.ncsu.edu/\~\,aszanto}

\date{\today}


\begin{abstract}
In 1853 Sylvester introduced
a {\em family\/} of double sum expressions
for two finite sets of indeterminates
and showed that some members of the family are essentially
the polynomial subresultants
of the monic polynomials associated with these sets.
A question naturally arises:
What are the {\em other\/} members of the family?
This paper provides a complete answer to this question.
The technique that we developed to answer the question
turns out to be general enough to characterise {\em all\/} members
of the family, providing a uniform method.
\end{abstract}

\end{frontmatter}

\section{Introduction}
 Let $A$ and
$B$ be finite lists (ordered sets) of distinct indeterminates. In
\cite{sylv}, Sylvester introduced for each $0\le p\le |A|, 0\le
q\le |B|$ the following {\em double-sum} expression in $A$ and
$B$:

$$
\operatorname{Sylv}^{p,q}(A,B;x)  :=
\sum_{
\begin{array}{l}
A^{\prime }\subset A,\,B^{\prime}\subset
B
\\[-1mm]
|A^{\prime}|=p,\,|B^{\prime}|=q
\end{array}}R(x,A^{\prime
})\,R(x,B^{\prime})\,\frac{R(A^{\prime},B^{\prime})\,R(A\backslash
A^{\prime},B\backslash B^{\prime})}{R(A^{\prime},A\backslash
A^{\prime })\,R(B^{\prime},B\backslash B^{\prime})},$$ where
\[
R(Y,Z):=\prod_{y\in Y,z\in Z}(y-z), \ \quad R(y,Z):=\prod_{z\in
Z}(y-Z).
\]
Let now  $f, g$ be univariate polynomials such that
$$\begin{array}{ccccc}
f&:=&\prod_{\alpha\in A}(x-\alpha)&=&x^m+a_{m-1}x^{m-1}+\ldots+a_0\\
g&:=&\prod_{\beta\in B}(x-\beta)&=&x^n+b_{n-1}x^{n-1}+\ldots+b_0,
\end{array}$$
where $m:=|A|\geq1$ and $n:=|B|\geq1$.


Since the double sum expressions are polynomials in $x$ and
symmetric in the $\alpha$'s and $\beta$'s, they can be expressed as
polynomials in $x$ whose coefficients are rational functions in the
$a_i$'s and the $b_j$'s.
In \cite{sylv}, the rational expression for
$\operatorname{Sylv}^{p,q}(A,B;x)$ is determined for the following
values of $(p,q)$, where without loss of of generality we assume $
1\le m\le n$ and we set $d:=p+q$ (see also \cite{LP}):
\begin{enumerate}
\item[(1)] If $0\le d<m\le n,$ then
$$\operatorname{Sylv}^{p,q}(A,B;x)=(-1)^{p(m-d)}{d\choose p}
\Sres_d(f,g),$$ where  $\Sres_d(f,g)$ is the $d$-th subresultant of
the polynomials $f$ and $g$, whose definition is recalled in Formula
(\ref{defsub})  below (cf. \cite[Art.~21]{sylv} and also
\cite[Theorem~0.1]{LP}).

\item[(1')] If $m=d<n,$ then
$$\operatorname{Sylv}^{p,q}(A,B;x)={m\choose p}f(x)$$
(cf. \cite[Art.~21]{sylv}  and also \cite[Proposition~2.9~(i)]{LP}).
In fact, the $d$-th subresultant is also well defined  for $d=m<n$
 as $\Sres_m(f,g)=f$. This implies that Case (1') can
be seen as a special case of Case (1).

\item[(2)] If $m<d<n-1,$ then
$\operatorname{Sylv}^{p,q}(A,B;x)=0$ (cf.
\cite[Arts.~23~\&~24]{sylv}).

\item[(3)] If $m<d=n-1$, then
$\operatorname{Sylv}^{p,q}(A,B;x)$ is a ``numerical multiplier'' of
$f(x)$ (cf. \cite[Art. 25]{sylv}), but the ratio is not established.

\item[(4)] If $m=d=n,$ then
$$\operatorname{Sylv}^{p,q}(A,B;x)={m-1\choose q}f(x)+{m-1\choose p}g(x)$$
(cf.\cite[Art.~22]{sylv} and also
\cite[Proposition~2.9~(ii)]{LP}).

\end{enumerate}

\bigskip
\noindent This note provides two contributions to this subject:

\medskip
\begin{itemize}
\item {\bf Completion}: Note that the above cases do not completely cover all possible values of
$p$ and $q$ such that $0\le p\le m$ and $0\le q\le n$:  values when
$n\leq p+q \leq m+n$  (except if $p+q=m=n$) are not covered. In Main
Theorem below we provide expressions for {\em all\/} the possible
values of $p$ and $q$, finally completing the previous efforts.

\item {\bf Uniformity}: \cite{sylv} and  \cite{LP} gave different proofs for each of the cases
listed above.
In Section~\ref{proof}, we provide a {\em uniform} technique
that can be applied to {\em all} the possible cases.
We obtained this technique by generalizing
the matrix formulation, used in \cite{DHKS} for dealing with the cases (1) and (2),
into a ``global'' matrix formulation. Approaches  to double sum expressions via matrix constructions have already been used
in \cite{Bch60,Bch78} (see also \cite{AJ06}).
\end{itemize}

\medskip
In order to state our main result we recall that for $0\le k\le m<
n$ or $0\le k < m=n$, the $k$-th subresultant of the polynomials $f$
and $g$ is defined as

\begin{equation}\label{defsub}
\Sres_k(f,g)
:=\det%
\begin{array}{|cccccc|c}
\multicolumn{6}{c}{\scriptstyle{m+n-2k}}\\
\cline{1-6}
a_{m} & \cdots & & \cdots & a_{k+1-\left(n-k-1\right)}& x^{n-k-1}f(x)&\\
& \ddots & && \vdots  & \vdots &\scriptstyle{n-k}\\
&  &a_{m}&\cdots &a_{k+1}& x^0f(x)& \\
\cline{1-6}
b_{n} &\cdots & & \cdots & b_{k+1-(m-k-1)}&x^{m-k-1}g(x)&\\
&\ddots &&&\vdots & \vdots &\scriptstyle{m-k}\\
&& b_{n} &\cdots & b_{k+1} & x^0g(x)&\\
\cline{1-6} \multicolumn{2}{c}{}
\end{array}
\end{equation}
with $a_{\ell}=b_{\ell}=0$ for $\ell<0$.

\smallskip
Expanding the determinant by the last column gives an
expression
\begin{equation}\label{subres}
\Sres_k(f,g)=F_k (x) f(x)+G_k(x)g(x)
\end{equation}
where the polynomials $F_k$ and $G_k$ (satisfying $\deg F_k\le
n-k-1$ and $\deg G_k\le  m-k-1$) are given by:

{\scriptsize
$$
\begin{array}{lclcrcl}  F_k &
:= &\det%
\begin{array}{|cccccc|}
\cline{1-6}
a_{m} & \cdots & & \cdots & a_{k+1-\left(n-k-1\right)}& x^{n-k-1}\\
& \ddots & && \vdots  & \vdots \\
&  &a_{m}&\cdots &a_{k+1}& x^0 \\
\cline{1-6}
b_{n} &\cdots & & \cdots & b_{k+1-(m-k-1)}&0\\
&\ddots &&&\vdots & \vdots \\
&& b_{n} &\cdots & b_{k+1} & 0\\
\cline{1-6}
\end{array}&,&
G_k
&:=&\det%
\begin{array}{|cccccc|}
\cline{1-6}
a_{m} & \cdots & & \cdots & a_{k+1-\left(n-k-1\right)}& 0\\
& \ddots & && \vdots  & \vdots \\
&  &a_{m}&\cdots &a_{k+1}& 0 \\
\cline{1-6}
b_{n} &\cdots & & \cdots & b_{k+1-(m-k-1)}&x^{m-k-1}\\
&\ddots &&&\vdots & \vdots \\
&& b_{n} &\cdots & b_{k+1} & x^0\\
\cline{1-6}
\end{array}\end{array}.
$$}

Now we are ready to state the main result that will be proven in the
next section.

\begin{undef}{Main Theorem.\ } \label{sylvester} Let $1\le m\le n,
\;\; 0\le p\le m, \;\;0\le q\le n,\;$  and set $d:=p+q,\  k :=
m+n-d-1,\ \sigma  := q(m-p)+n(d-m)+d+n-q-1 $. Then
 {
$$\operatorname{Sylv}^{p,q}(A,B;x)=\left\{
\begin{array}{lll}
(-1)^{p(m-d)}{d\choose p}
\Sres_{d}(f,g)&\mbox{ \ for \ }& 0\le d<  m \;\mbox{ or } \;m=d<n\\[1mm] 0&\mbox{ \ for \ }& m<d<n-1\\[1mm]
(-1)^{(m+q)(p+1)}{m\choose p} f&\mbox{ \ for \ }& m<d=n-1\\[1mm]
(-1)^{\sigma} \Big( {k\choose m-p}F_k \ f- {k\choose
n-q}G_k \ g\Big)&\mbox{ \ for \ }&n\le d\le  m+n-1\\[1mm]
\Res(f,g)\ fg&\mbox{ \ for \ } &d=m+n.
\end{array}
\right. $$}
\end{undef}

\medskip
We note that the previous Case (4) is covered here by the case $n\le
d\le m+n-1$, where indeed,  for $m=d=n$, $F_{m-1}=-1$ and
$G_{m-1}=1$.

\medskip
Finally, let us add two remarks kindly pointed out by one of the
referees. First, there are some avenues for further extensions: such
as more general types of summations considered by Sylvester himself
(see for instance \cite{sylv2}), and the cases when $A\cap
B\neq\emptyset$, or when $A$ and/or $B$ has repeated elements.
Second, the results concerning Sylvester's sums can be viewed as
generalizations of interpolation formulas (the $m=n-1$, $p=0$, $q=m$
case giving the Lagrange interpolation formula), and thus an
approach using specialization instead of linear algebra may also be
possible for proving such equalities.

\section{Proof}\label{proof}

As in \cite{DHKS}, we define for a polynomial $h(t)$, a finite
list
 $\Gamma:=(\gamma_1,\ldots,\gamma_u)$ of scalars  and
a non-negative integer $v$   the (non necessarily square) matrix
of size $v\times u$:
$$\langle h(t),\Gamma \rangle_v:=\begin{array}{|ccc|c}
\multicolumn{3}{c}{\scriptstyle{u}}&\\
\cline{1-3}
\gamma_1^0h(\gamma_1)&\dots & \gamma_u^0h(\gamma_u)&\\
\vdots& & \vdots& \scriptstyle v\\
\gamma_1^{v-1}h(\gamma_1)&\dots & \gamma_u^{v-1}h(\gamma_u)& \\
 \cline{1-3}
 \multicolumn{2}{c}{}
\end{array}.$$
For instance, under  this notation,  $$\langle
x-t,\Gamma\rangle_v= (\gamma_j^{i-1}x-\gamma_j^i)_{1\le i\le
v,1\le j\le u},$$ and  for $v=u$ we have the following equality
for the Vandermonde determinant $\mathcal{V}(\Gamma)$ associated
to $\Gamma$:
$$\mathcal{V}(\Gamma):=\det\big(\gamma_j^{i-1}\big)_{1\le i,j\le u}=\det\big(\langle 1,
\Gamma\rangle_u\big).$$

\par
\medskip
For the rest of the paper, $d\in\N$, $0\leq d\leq m+n$ and
$d':=m+n-d$. We take a new variable $T$ and   we denote by
$U_d(x,T)$  the following square matrix of size $m+n=d'+d$:
$$
U_d(x,T):=
\begin{array}{|c|c|l}
\multicolumn{1}{c}{\scriptstyle n}&\multicolumn{1}{c}{\scriptstyle m}&\\
\cline{1-2}
\langle 1, B\rangle_{d'}&\langle T, A\rangle_{d'}&\scriptstyle{d'}\\
\cline{1-2} \langle x-t,B\rangle_d
& \langle x-t, A\rangle_d&\scriptstyle d\\
\cline{1-2} \multicolumn{1}{c}{}
\end{array},
$$
where $\langle T, A\rangle_{d'}=\left(T\alpha^j\right)_{\alpha\in
A,\,0\leq j\leq d'-1}$. Finally we denote  by $u_d(x,T)$ its
determinant, that we develop in the powers of $T$:
\begin{equation}\label{udp}
u_d(x,T):=\det\big(U_d(x,T)\big)=u_{d,0}(x)T^m+ \cdots + u_{d,m-1}(x)T +
u_{d,m}(x).
\end{equation}

We are now ready to state our first result, that relates
$\operatorname{Sylv}^{p,d-p}(A,B;x)$ to the coefficient
$u_{d,p}(x)$:

\begin{thm}\label{matrix}
Let $0 \leq d \leq m+n $, $0\le p \le m$ and define $q:=d-p$.
Following Notation~(\ref{udp}), we have that if $0\leq q \leq n$
then
$$
u_{d,p}(x)=(-1)^{q(m-p)}\V(A)\,\V(B)\,\operatorname{Sylv}^{p,q}(A,B;x)
$$
while otherwise  $u_{d,p}(x)=0$.
\end{thm}

\medskip
\begin{pf}
For any set $A'\subset A$ (resp. $B'\subset B$) we will denote with $A''$ (resp. $B''$) its complementary set, i.e.
$A'':=A\setminus A'$ (resp. $B'':=B\setminus B'$).
\par
We perform  a
Laplace expansion of the determinant of the matrix $U_d(x,T)$ on
the
 last $d$ rows and we get the following
expression:
$$
u_d(x,T)=\sum_{\begin{array}{c}
A'\subset A,B'\subset B\\[-1mm]
|A'|+ |B'|=d
\end{array}}\sigma(B''\cup A'',B\cup A)T^{m-|A'|}\V(B''\cup A'') R(x,B')R(x,A')\V(B'\cup A'),$$
 where, as
in \cite{DHKS}, ``$\cup$" stands for list concatenation,
``$\backslash$" means list subtraction and, for $S\subseteq T$
finite lists, $\sigma\left( S,T\right) :=\left( -1\right) ^{j},\
j$ being the number of transpositions needed to take $T$ to $S
\cup (T\backslash S)$.
\par
We write $u_d(x,T)$ in powers of $T$, with $0\le p\le m$ and $0\le
q:=d-p\le n$ implying $\max\{0,d-n\}\le p\le \min\{d,m\}$:
{\small $$
u_d(x,T)=\sum_{p=\max\{0,d-n\}}^{\min\{d,m\}}\left(\sum_{{\tiny\begin{array}{c}
A'\subset A,B'\subset B\\[-1mm]
|A'|=p,\,|B'|=q
\end{array}}}\sigma(B''\cup A'',B\cup A)R(x,B')R(x,A')
\V(B'\cup A')\V(B''\cup A'')\right)T^{m-p}.
$$}
We recall the elementary fact that transposing  a block of $j$
columns with an adjacent  block of $i$ columns produces in the
determinant a change of sign of order $(-1)^{ij}$. Hence, for
$|A'|=p$ and $|B'|=q$, {\small $$\sigma (B''\cup
A'',B\cup A)=\sigma(A'',A)\sigma(B'',B)(-1)^{q(m-p)},$$} and we have for $\max\{0,d-n\}\le p\le
\min\{d,m\}$:
$$
u_{d,p}(x)=(-1)^{q(m-p)}\sum_{\begin{array}{c}
A'\subset A,B'\subset B\\[-1mm]
|A'|=p,\,|B'|=q
\end{array}}\sigma(A'',A)\sigma(B'',B)R(x,A')
R(x, B')\V(B'\cup A')\V(B''\cup A'').$$ Now
we apply repeatedly  the elementary fact that
$$
  \V(X\cup
Y)=\V(X)\,\V(Y)\,R(Y,X)$$ for any pair of finite lists $X,Y$:

{ \begin{eqnarray*} \V(B'\cup A')\,\V(B''\cup
A'')&=&\V(A')\V(B')R(A',B')\V(A'')\V(B'')R(A'',B'')\\
&=&\frac{\V(A''\cup A')}{R(A',A'')}\frac{\V(B''\cup B')}{R(B',B'')}R(A',B')R(A'',B'')\\
&=& \sigma(A'',A)\sigma(B'',B)\V(A)\V(B)\frac{R(A',B')R(A'',B'')}{R(A',
A'')R(B', B'')}.
\end{eqnarray*}}
We finally obtain that {
\begin{eqnarray*} u_{d,p}(x)&=&
(-1)^{q(m-p)}\V(A)\V(B)\sum_{\begin{array}{c}
A'\subset A,B'\subset B\\[-1mm]
|A'|=p,\,|B'|=q
\end{array}}R(x,A')
R(x, B')\frac{R(A',B')R(A'',B'')}{R(A',A'')R(B',B'')}\\
&=&
(-1)^{q(m-p)}\V(A)\V(B)\operatorname{Sylv}^{p,q}(A,B;x).
\end{eqnarray*} }
\end{pf}

In view of Theorem~\ref{matrix}, in order to produce a rational
expression for $\operatorname{Sylv}^{p,q}(A,B;x)$ it is
enough to give a rational expression for $u_{d,p}(x)$.  To this
aim we first observe the following straightforward  factorization
formula for $U_d(x,T)$ as a product of two rectangular matrices of
sizes $(m+n)\times (m+n+1)$ and $(m+n+1)\times (m+n)$
respectively.

\begin{lem} \
{\small
\begin{equation}\label{factor1}
 U_d(x,T)=
\begin{array}{r|c|c|}
\multicolumn{1}{c}{}&\multicolumn{1}{c}{\scriptstyle{d'}}&\multicolumn{1}{c}{\scriptstyle
d+1}\\
\cline{2-3} \scriptstyle d'& {\II}_{d'}
&{\bf 0}\\
\cline{2-3} \scriptstyle{d}&{\bf 0}&\begin{array}{rrrrrr}
x&-1&0&\ldots&0&0\\
0&x&-1&\ldots&0&0\\
&&\ddots&\ddots&&\\
&&&\ddots&\ddots&\\
0&0&0&\ldots&x&-1
\end{array}\\
\cline{2-3} \multicolumn{1}{c}{}
\end{array}
\ \
\begin{array}{|c|c|l}
\multicolumn{1}{c}{\scriptstyle n}&\multicolumn{1}{c}{\scriptstyle m}&\\
\cline{1-2} \langle 1, B\rangle_{d'}
&\langle T, A\rangle_{d'}&\scriptstyle{d'}\\
\cline{1-2} \langle
1,B\rangle_{d+1}& \langle 1, A\rangle_{d+1}&\scriptstyle{d+1}\\
\cline{1-2} \multicolumn{1}{c}{}
\end{array},
\end{equation} }
where $\II_{d'}$ denotes the identity matrix of size $d'$.
\end{lem}

\smallskip

The previous factorization of $U_{d}(x,T)$ immediately yields

\smallskip
\begin{prop}[Arts.~23~\&~24 \cite{sylv}]\label{cero} \
 If $m<d<n-1$,  then $u_d(x,T)=0.$
\end{prop}
\begin{pf}
The assumption implies $\max\{d',d+1\}<n$. Then the first $n$
columns of the matrix at the right of (\ref{factor1}) have
deficient rank since all $n\times n$ minors vanish.  A
Binet-Cauchy expansion of $u_d(x,T)$ therefore implies that
$u_d(x,T)$ vanishes as well.
\end{pf}

\smallskip
Our goal now is to provide a factorization  like in
(\ref{factor1}), but with square matrices, that allows to recover
$u_d(x,T)$.

\smallskip

 \begin{thm}\label{factor2} Let $1\le m\le n$.
If $0\le d\le m$ or $n-1\le d\le  m+n$, then  there exist
polynomials
 $P(x):=P_0+\cdots +P_dx^d$ and
$Q(x):=Q_0+\dots +Q_{d'-1}x^{d'-1}$ with $P\ne 0$, $k:=\deg P\le d$
and  $\deg Q\le d'-1$ if $d'\ne 0$ such that  we have the following
 matrix identity:
{\scriptsize $$
\begin{array}{r|c|c|}
\multicolumn{1}{c}{}&\multicolumn{1}{c}{\scriptstyle{d'}}&\multicolumn{1}{c}{\scriptstyle
d+1}\\
\cline{2-3} \scriptstyle d'&\begin{array}{ccc}& {\II_{d'}}
&\end{array}
&{\bf 0}\\
\cline{2-3} \scriptstyle{d}&\bf{0}&\begin{array}{cccc}
x&-1& &0\\
&\ddots&\ddots& \\
0&\dots & x&-1
\end{array}\\
\cline{2-3} \scriptstyle{1}&\begin{array}{ccc} Q_0&\dots &Q_{d'-1}
\end{array}&\begin{array}{lccr}
P_0&\dots &\dots  &P_d
\end{array}\\
\cline{2-3}\multicolumn{1}{c}{}
\end{array}
\
\begin{array}{|c|c|c|l}
\multicolumn{1}{c}{\scriptstyle n}&\multicolumn{1}{c}{\scriptstyle
m}&
\multicolumn{1}{c}{\scriptstyle 1}&\\
\cline{1-3} \langle 1, B\rangle_{d'}
&\langle T, A\rangle_{d'}&\bf 0& \scriptstyle{d'}\\
\cline{1-3} \langle
1,B\rangle_{d+1}& \langle 1, A\rangle_{d+1}&{\bf e}_k &\scriptstyle{d+1}\\
\cline{1-3} \multicolumn{1}{c}{}
\end{array} =
 \begin{array}{|c|c|l}
\multicolumn{1}{c}{\scriptstyle{m+n}}&\multicolumn{1}{c}{\scriptstyle
1}&\multicolumn{1}{c}{}\\
\cline{1-2} U_d(x,T)
&{\bf *}& \scriptstyle m+n\\
\cline{1-2} \bf 0 & P_k&\scriptstyle{1}\\
\cline{1-2} \multicolumn{1}{c}{}
\end{array} ,
$$ }
\noindent where ${\bf e}_k$ is defined  as the vertical vector
 of size $d+1$ with a single non-zero entry $1$ in position $k +1$
 and $P_k$ is the leading coefficient of $P$.

Moreover, $P(x)$ can be defined as
 {$$P(x)=\left\{
\begin{array}{lcl}\Sres_d(f,g)&\mbox{ \ for \ }& 0\le d< m \;\;\mbox{ or
}\;\; d=m<n
\\  f  &\mbox{ \ for \ } & m<d=n-1\\  F_{d'-1}f+TG_{d'-1}g
&\mbox{for}& n\le d< m+n
\\  fg  &\mbox{ \ for \ } & d=m+n,\end{array}\right.
 $$}
where $F_{d'-1},G_{d'-1}$  are as in Identity~(\ref{subres}) for
$k=d'-1$.
\end{thm}

\begin{rem} We note that
$P(x)=F_{d'-1} f+TG_{d'-1}g$ is  the determinant of a matrix
similar to the matrix (\ref{defsub}) that defines
$\Sres_{d'}(f,g)$: we simply need to replace $g(x)$ by $Tg(x)$ in
the last column of the matrix (\ref{defsub}).
\end{rem}

\begin{pf} To get the  factorization stated in Theorem \ref{factor2},  we only
need to look at
the equations that can be read from the lower row of
the matrix on the right. These are
$$\begin{array}{lcl} ({Q_0}+\dots +{Q_{d'-1}}\beta^{d'-1}) + ({P_0}
+\dots +
P_{d} \beta^d )&=& 0 \\
( T {Q_0} +\dots + T {Q_{d'-1}} \alpha^{d'-1}) + ({P_0}  +\dots +
{P_d}\alpha^{d})&=&0\end{array} $$
for all $\beta \in B$, $\alpha
\in A$. In order to solve these equations, it is enough to produce polynomials
$P(x):=P_0+\cdots +P_dx^d$ and $Q(x):=Q_0+\cdots +Q_{d'-1}x^{d'-1}$
with $P\ne 0$, $\deg P\le d$ and $\deg Q\le d'-1$ if $d'\ne 0$ such
that the following $m+n$ equations are satisfied:
\begin{equation}\label{cond}
\left\{ \begin{array}{ll} Q(\beta)+ P(\beta)=0, & \forall \beta
\in B\\TQ(\alpha)+P(\alpha)=0,&\forall \alpha \in A.
\end{array}\right.
\end{equation}

\smallskip

For ${0\le d\le m}$ if $m<n$ and $0\le d<m$ if $m=n$, we define $$
\left\{ \begin{array}{lcl} P(x)&:=& \Sres_d(f,g)\ =\
F_d(x)f(x)+G_d(x)g(x)\\[2mm]
Q(x)&:=&-F_d(x)f(x) -\frac{1}{T}G_d(x)g(x)\end{array}\right.$$
where $F_d,G_d$ are as in Identity~(\ref{subres}) for $k:=d$. Thus
$\deg P=\deg \Sres_d(f,g)=d$ and $\deg_x Q\le \max\{\deg
(F_df),\deg (G_dg)\}\le  d'-1$. We look at Condition~(\ref{cond}):
$$\left\{ \begin{array}{lcl}
 Q(\beta)+ P(\beta)&=& (1-\frac{1}{T})G_d(\beta)g(\beta )\quad  = \quad  0,
  \quad \forall \beta \in B\\[1mm]
TQ(\alpha)+ P(\alpha)&=&(1-T)F_d(\alpha)f(\alpha) \quad = \quad 0,
\quad \forall \alpha \in A.
\end{array}\right.$$

\smallskip
For $m<d=n-1$, we define $$P(x):=f(x) \quad \mbox{and} \quad
Q(x):=-f(x).$$ We have $\deg P=m< d$ and $\deg Q=m= m+n-d-1=d'-1$
in this case. Condition~(\ref{cond}) is trivially satisfied.

\smallskip
For ${n\le d< m+n}$, we observe that $0\le d'-1\le m-1$. Thus
$\Sres_{d'-1}(f,g)$ is well defined and we define
 $$ \left\{ \begin{array}{lcl} Q(x)&:=& -\Sres_{d'-1}(f,g)\ =
-F_{d'-1}(x)f(x)-G_{d'-1}(x)g(x)\\[2mm]
P(x)&:=&F_{d'-1}(x)f(x) + TG_{d'-1}(x)g(x)\end{array}\right.$$
where $F_{d'-1},G_{d'-1}$ are as in Identity~(\ref{subres}) for
$k:=d'-1$. Thus $\deg Q=\deg \Sres_{d'-1}(f,g)=d'-1$ and   $\deg_x
P\le \max\{\deg(F_{d'-1}f),\deg(G_{d'-1}g)\}\le
\max\{m+n-(d'-1)-1,n+m-(d'-1)-1\}=d$. Also $P\ne 0$  since the
leading terms can not cancel each other. We look again at
Condition~(\ref{cond}):
$$\left\{\begin{array}{lcl}
 Q(\beta)+ P(\beta)&=&  (T-1)G_{d'-1}(\beta)g(\beta )\quad  = \quad  0,
  \quad \forall \beta \in B\\[1mm]
TQ(\alpha)+ P(\alpha)&=&(1-T)F_{d'-1}(\alpha)f(\alpha)\quad =\quad
0,\quad \forall \alpha \in A.
\end{array}\right.$$

\smallskip
For $d=m+n$, since $d'=0$ in this case, we define $P(x)=f(x)g(x)$,
which is  of degree $d$,  to satisfy Condition~(\ref{cond}).
\end{pf}

Theorem~\ref{factor2} immediately implies that $u_d(x,T)$ can be
computed as the determinant of two square matrices for  the values
of $d\le m$ and $ n-1\le d$. Our  next goal is to
 compute $P_k$ in each case, as well as
 the determinants of these square matrices.\\
To this aim,  for $0\le d\le m<n$ or $0\le d<m=n$, we set
$\Delta_k(f,g)$ for the leading coefficient of
 $\Sres_k(f,g)$, i.e. $\Delta_k(f,g)$ is
the $k$-th scalar subresultant of $f,g$.

For $k=m=n$, we define for the coherence of the next results
$$\Delta_m(f,g):=1.$$

\begin{lem}\label{peka}
Let $1\le m\le n$. Following the notations of Theorem~\ref{factor2},
we have {\small
$$\left\{ \begin{array}{lclcl}
\deg P=d  &\mbox{and}& P_k=\Delta_d(f,g)&  \mbox{ \ for \ }& 0\le
d\le m<n \mbox{ \ or \  }
d<m=n\\
\deg P=m &\mbox{and}& P_k=1&  \mbox{ \ for \ }& m<d=n-1\\
\deg P=d  &\mbox{and} & P_k=(-1)^{d-n}\Delta_{d'}(f,g)(T-1) & \mbox{
\ for \ }&m\le n\le d< m+n
\\
\deg P=d &\mbox{and}& P_k=1&  \mbox{ \ for \ }& d=m+n.
\end{array}\right.$$}
\end{lem}
\begin{pf} The first two cases and the last case are straightforward from the
definition of $P_k$.\\
 For $ m\le n\le d< m+n$,
we have that $P(x)=F_{d'-1}(x)f(x)+TG_{d'-1}(x)g(x)$. Thus
 $\deg_x P=
\max\{\deg_x(F_{d'-1}f),\deg_x(G_{d'-1}g)\}$ since the leading
terms can not cancel each other. A direct computation on the
matrix in (\ref{defsub}) that defines  $\Sres_k(f,g)$ shows that
---since for $k:=d'-1<m$, $n-k>1$ and $m-k>1$ hold--- then
$\deg_xF_k=n-k-1=d-m$
and $\deg_xG_k=m-k-1=d-n$. Therefore $\deg_xP=\max\{m+n-(d'-1)-1,n+m-(d'-1)-1\}=d$.\\
Finally,   since $f$ and $g$ are monic, the leading coefficient of
$F_k(x)$ equals
$$(-1)^{m+n-2k+1}(-1)^{n-k+1}\Delta_{k+1}(f,g)=
(-1)^{d-n+1}\Delta_{d'}(f,g)$$ and the leading coefficient of
$G_k(x)$ equals
$$(-1)^{m+n-2k+n-k+1}\Delta_{k+1}(f,g)=
(-1)^{d-n}\Delta_{d'}(f,g).$$ Therefore $P_d=
(-1)^{d-n}\Delta_{d'}(f,g)(T-1)$.
\end{pf}

\begin{lem}\label{companion}
{}

$$\det\begin{array}{|c|c|l}
\multicolumn{1}{c}{\scriptstyle{d'}}&\multicolumn{1}{c}{\scriptstyle
d+1}&\multicolumn{1}{c}{}\\
\cline{1-2} \begin{array}{ccc}& {\II_{d'}} &\end{array} &{\bf
0}&\scriptstyle
d'\\
\cline{1-2} \bf{0}&\begin{array}{cccc}
x&-1& &0\\
&\ddots&\ddots& \\
0&\dots & x&-1
\end{array}&\scriptstyle{d}\\
\cline{1-2} \begin{array}{ccc} Q_0&\dots &Q_{d'-1}
\end{array}&\begin{array}{lccr}
P_0&\dots &\dots  &P_d
\end{array}&\scriptstyle{1}\\
\cline{1-2}\multicolumn{1}{c}{}
\end{array}=P_0+\cdots + P_dx^d=P(x).
$$
\end{lem}
\begin{pf}Because of the block triangular structure,
this determinant equals $$\det \begin{array}{|c|l}
\multicolumn{1}{c}{\scriptstyle d+1}&\multicolumn{1}{c}{}\\
\cline{1-1} \begin{array}{cccc}
x&-1& &0\\
&\ddots&\ddots& \\
0&\dots & x&-1
\end{array}&\scriptstyle{d}\\
\cline{1-1} \begin{array}{lccr} P_0&\dots &\dots  &P_d
\end{array}& \scriptstyle{1}\\
\cline{1-1}\multicolumn{1}{c}{}
\end{array}.$$ We can permute the first $d$-block
with the last row and expand the determinant by this new first
row. We get
$$(-1)^d\left(P_0(-1)^d- P_1 x(-1)^d + \dots +(-1)^{d} P_d x^d\right).$$
\end{pf}

\begin{lem}\label{ion} Let $1\le m\le n$. Then
{
\begin{eqnarray*}&&\det \begin{array}{|c|c|c|l} \multicolumn{1}{c}{\scriptstyle n}&\multicolumn{1}{c}{\scriptstyle m}&
\multicolumn{1}{c}{\scriptstyle 1}&\\
\cline{1-3} \langle 1, B\rangle_{d'}
&\langle T, A\rangle_{d'}&\bf 0& \scriptstyle{d'}\\
\cline{1-3} \langle
1,B\rangle_{d+1}& \langle 1, A\rangle_{d+1}&{\bf e}_d &\scriptstyle{d+1}\\
\cline{1-3} \multicolumn{1}{c}{}
\end{array}=\\&&=\left\{
\begin{array}{lcl}
(-1)^{dm}\V(A)\,\V(B)\,\Delta_d(f,g)\,T^{m-d}(T-1)^d&\mbox{ \ for \
}& 0\le d\le
m\\[1mm]
(-1)^{m(d-1)+d}\V(A)\V(B)(T-1)^m&\mbox{ \ for \ }& m<d=n-1\\[1mm]
(-1)^{d'n}\V(A)\,\V(B)\,\Delta_{d'}(f,g)\,(T-1)^{d'}&\mbox{ \ for \
}&n\le d< m+n
\\[1mm]
\V(A)\,\V(B)\,\Res(f,g)&\mbox{ \ for \ }&d= m+n.
\end{array}
\right. \end{eqnarray*}}
\end{lem}

\begin{pf}
First, let us recall \cite[Lemma~2]{DHKS}:
$$
\Sres_k(f,g){\mathcal V}(A)=\det\begin{array}{|c|c}
\multicolumn{1}{c}{ \scriptstyle{m}}&\\[1mm]
\cline{1-1}
\langle x-t,A\rangle_k&\scriptstyle{k}\\[1mm]
\cline{1-1}
\langle g(t),A\rangle_{m-k}&\scriptstyle{m-k}\\[1mm]
\cline{1-1} \multicolumn{2}{c}{}
\end{array},
$$
which  implies that its leading coefficient satisfies
\begin{equation}\label{sigue}
\Delta_k(f,g){\mathcal V}(A)=\det\begin{array}{|c|c}
\multicolumn{1}{c}{ \scriptstyle{m}}&\\[1mm]
\cline{1-1}
\langle 1,A\rangle_k&\scriptstyle{k}\\[1mm]
\cline{1-1}
\langle g(t),A\rangle_{m-k}&\scriptstyle{m-k}\\[1mm]
\cline{1-1} \multicolumn{2}{c}{}
\end{array}.
\end{equation}
To simplify the notation of the proof, we will denote  the matrix on
the left side of the claim of the Lemma by $M_d$.\\
 In case $0\le d\le m$ or $n\le d< m+n$, $\deg
P(x)=d$ by Lemma \ref{peka} and ${\bf e}_d:=(0, \dots ,0 ,1)^t$.
Therefore
$$\det(M_d)= \det \begin{array}{|c|c|l}
\multicolumn{1}{c}{\scriptstyle n}&\multicolumn{1}{c}{\scriptstyle m}&\\
\cline{1-2} \langle 1, B\rangle_{d'}
&\langle T, A\rangle_{d'}&\scriptstyle{d'}\\
\cline{1-2} \langle
1,B\rangle_{d}& \langle 1, A\rangle_{d}&\scriptstyle{d}\\
\cline{1-2} \multicolumn{1}{c}{}
\end{array}.
$$
For $d\le m$, we have that  $d'\ge n\ge m\ge d$ holds and
therefore row operations yield {
\begin{eqnarray*} \det(M_d)&=& \det \begin{array}{|c|c|l}
\multicolumn{1}{c}{\scriptstyle n}&\multicolumn{1}{c}{\scriptstyle m}&\\
\cline{1-2} \langle 1, B\rangle_{d'}
&\langle T, A\rangle_{d'}&\scriptstyle{d'}\\
\cline{1-2} \bf 0 & \langle 1-T, A\rangle_{d}&\scriptstyle{d}\\
\cline{1-2} \multicolumn{1}{c}{}
\end{array}
 \\
 &=&\det  \begin{array}{|c|c|l}
\multicolumn{1}{c}{\scriptstyle n}&\multicolumn{1}{c}{\scriptstyle m}&\\
\cline{1-2} \langle 1, B\rangle_{n}
&\langle T, A\rangle_{n}&\scriptstyle{n}\\
\cline{1-2} \bf 0
&\langle T g(t), A\rangle_{m-d}&\scriptstyle{m-d}\\
\cline{1-2} \bf 0 & \langle 1-T, A\rangle_{d}&\scriptstyle{d}\\
\cline{1-2} \multicolumn{1}{c}{}
\end{array}\quad \mbox{since for all } \beta\in B,
g(\beta)=0 \\
&=& \V(B)T^{m-d}(1-T)^d \det
\begin{array}{|c|l}
\multicolumn{1}{c}{\scriptstyle m}&\\
\cline{1-1}
\langle g(t),A\rangle_{m-d}&\scriptstyle{m-d}\\
\cline{1-1}
\langle 1, A\rangle_{d}&\scriptstyle{d}\\
\cline{1-1} \multicolumn{1}{c}{}
\end{array}
\\ &=& \V(B)\,T^{m-d}(1-T)^d(-1)^{d(m-d)}\,\V(A)\,\Delta_d(f,g) \qquad \mbox{by
(\ref{sigue})}\\
&=& (-1)^{d(m-d+1)}\V(A)\,\V(B)\,\Delta_d(f,g)\,T^{m-d}(T-1)^d\\
&=&(-1)^{dm}\V(A)\,\V(B)\,\Delta_d(f,g)\,T^{m-d}(T-1)^d.
\end{eqnarray*}}
In case $d\ge n$, we have that $d'\le m\le d$ holds and therefore
row operations yield

{
\begin{eqnarray*} \det(M_d)&=& \det
\begin{array}{|c|c|l}
\multicolumn{1}{c}{\scriptstyle n}&\multicolumn{1}{c}{\scriptstyle m}&\\
\cline{1-2}\bf 0
&\langle T-1, A\rangle_{d'}&\scriptstyle{d'}\\
\cline{1-2} \langle
1,B\rangle_{d}& \langle 1, A\rangle_{d}&\scriptstyle{d}\\
\cline{1-2} \multicolumn{1}{c}{}
\end{array}\\ &=& \det
\begin{array}{|c|c|l}
\multicolumn{1}{c}{\scriptstyle n}&\multicolumn{1}{c}{\scriptstyle m}&\\
\cline{1-2}\bf 0
&\langle T-1, A\rangle_{d'}&\scriptstyle{d'}\\
\cline{1-2} \langle
1,B\rangle_{n}& \langle 1, A\rangle_{n}&\scriptstyle{n}\\
\cline{1-2} \bf 0& \langle g(t), A\rangle_{d-n}&\scriptstyle{d-n}\\
\cline{1-2} \multicolumn{1}{c}{}
\end{array}\\&=&(-1)^{d'n}\V(B)\,(T-1)^{d'}\,\det
\begin{array}{|c|l}
\multicolumn{1}{c}{\scriptstyle n}&\\
\cline{1-1}
\langle 1, A\rangle_{d'}&\scriptstyle{d'}\\
\cline{1-1}
\langle g(t), A\rangle_{d-n}&\scriptstyle{d-n}\\
\cline{1-1} \multicolumn{1}{c}{}
\end{array}\\
&=&(-1)^{d'n}\V(A)\,\V(B)\,\Delta_{d'}(f,g)\,(T-1)^{d'}.
\end{eqnarray*}}

In case $m<d=n-1$, $\deg P=m$ and ${\bf e}_d$ is the vertical
vector with a single non-zero entry 1 in position $m+1$.  Since
$d+1=n$, $d'=m+1$ and $n\ge m+1$, \begin{eqnarray*} \det(M_d)&=& \det
\begin{array}{|c|c|c|l} \multicolumn{1}{c}{\scriptstyle
n}&\multicolumn{1}{c}{\scriptstyle m}&
\multicolumn{1}{c}{\scriptstyle 1}&\\
\cline{1-3} \langle 1, B\rangle_{m+1}
&\langle T, A\rangle_{m+1}&\bf 0& \scriptstyle{m+1}\\
\cline{1-3} \langle
1,B\rangle_{n}& \langle 1, A\rangle_{n}&{\bf e}_{d} &\scriptstyle{n}\\
\cline{1-3} \multicolumn{1}{c}{}
\end{array}\\&=&
\det
\begin{array}{|c|c|c|l} \multicolumn{1}{c}{\scriptstyle
n}&\multicolumn{1}{c}{\scriptstyle m}&
\multicolumn{1}{c}{\scriptstyle 1}&\\
\cline{1-3} {\bf 0}
&\langle T-1, A\rangle_{m+1}&-{\bf e}_d&\scriptstyle{m+1}\\
\cline{1-3} \langle
1,B\rangle_{n}& \langle 1, A\rangle_{n}&{\bf e}_{d}&\scriptstyle{n}\\
\cline{1-3} \multicolumn{1}{c}{}
\end{array}
\\&=& (-1)^{(m+1)n}
\det \langle 1,B\rangle_{n} \, \det
\begin{array}{|c|c|l} \multicolumn{1}{c}{\scriptstyle m}&
\multicolumn{1}{c}{\scriptstyle 1}&\\
\cline{1-2} \langle T-1, A\rangle_{m+1}&-{\bf e}_d&\scriptstyle{m+1}\\
\cline{1-2} \multicolumn{1}{c}{}
\end{array}\\&=& -(-1)^{(m+1)n}
\det \langle 1,B\rangle_{n} \,\det \langle 1,A\rangle_{m}(T-1)^m
\\[2mm]
&=& (-1)^{m(d-1)+d}\V(A)\V(B)(T-1)^m.
\end{eqnarray*}

Finally the case $d=m+n$ is straightforward since
$$\det(M_d)=\V(B\cup A)=\V(A)\V(B)\Res(f,g).$$
\end{pf}

We are ready now to compute $u_d(x,T)$ for all values of $d$,
$0\le d\le m+n$,  and to deduce
$\operatorname{Sylv}^{p,q}(A,B;x)$ for all possible
values of $p$ and $q$.

\begin{thm}\label{detUd} Let $1\le m\le n$. Then
{\scriptsize $$u_d(x,T)=\left\{
\begin{array}{lll}
(-1)^{dm}\V(A)\,\V(B)\,\Sres_d(f,g)\,T^{m-d}(T-1)^d&\mbox{ \ for \
}& 0\le d<  m  \mbox{ \ or \ }
m=d<n\\[1mm]
0&\mbox{ \ for \
}& m<d<n-1\\
(-1)^{\sigma}\V(A)\,\V(B)\,f(x)\,(T-1)^m&\mbox{ \ for \
}& m<d=n-1\\
(-1)^{\sigma}\V(A)\,\V(B)\,\left(F_{d'-1}(x)f(x) +
TG_{d'-1}(x)g(x)\right)\,(T-1)^{d'-1}&\mbox{ \ for \ }&n\le d< m+n
\\[1mm]
\V(A)\,\V(B)\,\Res(f,g)\,f(x)\,g(x)&\mbox{ \ for \ }& d=m+n,
\end{array}
\right.
$$}
where $\sigma=(d'-1)n+d$, and $F_{d'-1}, G_{d'-1}$ are defined as
in Identity~(\ref{subres}) for $k:=d'-1$.
\end{thm}

\begin{pf}  If  $m<d<n-1$ then by
Proposition~\ref{cero} we have that $u_d(x,T)=0$. For the other
cases of $0\le d\le m+n$, we apply Theorem~\ref{factor2} and Lemma
\ref{companion}. Using the notation of Theorem~\ref{factor2} we get
$$
 u_d(x,T)\cdot P_k=P(x) \cdot \det \begin{array}{|c|c|c|l}
\multicolumn{1}{c}{\scriptstyle n}&\multicolumn{1}{c}{\scriptstyle
m}&
\multicolumn{1}{c}{\scriptstyle 1}&\\
\cline{1-3} \langle 1, B\rangle_{d'}
&\langle T, A\rangle_{d'}&\bf 0& \scriptstyle{d'}\\
\cline{1-3} \langle
1,B\rangle_{d+1}& \langle 1, A\rangle_{d+1}&{\bf e}_k &\scriptstyle{d+1}\\
\cline{1-3} \multicolumn{1}{c}{}
\end{array}.
$$

Now for each of the following cases we also apply Lemmas
\ref{peka} and \ref{ion}:
\\
 For $0\le d< m$ or $d=m$ if
$m<n$,
 $P(x)=\Sres_d(f,g)$, $k=d$  and
$P_k=\Delta_d(f,g)$, therefore {
\begin{eqnarray*} u_d(x,T)&=&
\frac{1}{\Delta_d(f,g)}\left(\Sres_d(f,g)(-1)^{dm}\V(A)\,\V(B)\,\Delta_d(f,g)\,T^{m-d}(T-1)^d\right)
\\[1mm]
&=&
(-1)^{dm}\V(A)\,\V(B)\,\Sres_d(f,g)\,T^{m-d}(T-1)^d.\end{eqnarray*}}
\\For $m<d=n-1$ we have that $P(x)=f(x)$ $k=m$ and $P_k=1$, then
$$u_d(x,T)=f(x)(-1)^{m(d-1)+d}\V(A)\V(B)(T-1)^m ,$$
and to get the sign $(-1)^{\sigma}$ as in the claim, we note that
in this case $m=d'-1$ and thus $m(d-1)+d\equiv (d'-1)n+d \pmod 2$.

 For ${n\le d< m+n}$ we have that
$P(x)=F_{d'-1}(x)f(x)+TG_{d'-1}(x)g(x)$, $k=d'-1$ and $P_k=
(-1)^{d-n}\Delta_{d'}(f,g)(T-1)$. We conclude

{\small
\begin{eqnarray*} u_d(x,t)&=&
\frac{\left(F_{d'-1}(x)f(x) + T
G_{d'-1}(x)g(x)\right)(-1)^{nd'}\V(A)\,\V(B)\,\Delta_{d'}(f,g)\,(T-1)^{d'}}{(-1)^{d-n}(T-1)
\Delta_{d'}(f,g)}\\[1mm]
&=& (-1)^{n(d'-1)+d}\V(A)\,\V(B)\,\left(F_{d'-1}(x)f(x) +
TG_{d'-1}(x)g(x)\right)\,(T-1)^{d'-1}.\end{eqnarray*}}

The last case, $d=m+n$, is straightforward. We note that in this
case $u_x(x,T)$ is equal to $\V\big(A\cup B\cup\{x\}\big)$ up to a
sign.
\end{pf}

\begin{pf}{\bf (Main Theorem.)}\\
By Theorem \ref{matrix}
we have that
\begin{equation}\label{th2.1}
u_d(x,T)=(-1)^{q(m-p)}\V(A)\,\V(B)\,\operatorname{Sylv}^{p,q}(A,B;x).
\end{equation}
 For $0\le d:=p+q\le m<n$ or for $0\le d<m=n$, we have by
Theorem~\ref{detUd}
$$u_d(x,T)=\sum_{p=0}^{d}u_{d,p}(x)T^{m-p}=(-1)^{dm}\V(A)\,\V(B)\,\Sres_d(f,g)\,T^{m-d}(T-1)^d,$$
 which implies that \begin{eqnarray*}u_{d,p}(x)&=&(-1)^{dm}(-1)^p
{d\choose d-p}\,\V(A)\,\V(B)\,\Sres_d(f,g).
\end{eqnarray*}
Therefore, using (\ref{th2.1}),
\begin{eqnarray*}\operatorname{Sylv}^{p,q}(A,B;x)&=&(-1)^{dm+p-q(m-p)}{d\choose
p}\Sres_d(f,g)\\
&=&(-1)^{p(m-d)}{d\choose p}\Sres_d(f,g) \end{eqnarray*} since
{\small
$$dm+p-q(m-p)=pm +p+qp =p(m-d)+p(d+1+q)\equiv
p(m-d)+p(p+1)\pmod 2.$$} \\
For $m<d<n-1$, $\operatorname{Sylv}^{p,q}(A,B;x)=0$
since
$u_d(x,T)=0$.\\
For $m<d:=p+q=n-1$, $$u_d(x,T)=\sum_{p=0}^m
u_{d,p}(x)T^{m-p}=(-1)^{(d'-1)n+d}\V(A)\V(B)(T-1)^mf(x)$$ which
implies that
\begin{eqnarray*} u_{d,p}(x,T)&= &(-1)^{(d'-1)n+q}  {m\choose p}
\V(A)\V(B)f(x).
\end{eqnarray*}
Therefore, using (\ref{th2.1}), we get
$$\operatorname{Sylv}^{p,q}(A,B;x)=
(-1)^{(m+q)(p+1)}{m\choose p}f(x),$$
 since{\small
$$(d'-1)n+q-q(m-p)=m(p+q-1)+q-qm+qp\equiv (m+q)(p+1)\pmod 2.$$}
\\[2mm]
  For $m\le n\le d:=p+q<
m+n$, {\scriptsize
$$u_d(x,T)=\sum_{p=d-n}^{m}u_{d,p}(x)T^{m-p}=(-1)^{n(d-m)+d}\V(A)\,\V(B)\,\left(F_{d'-1}(x)f(x)
+ TG_{d'-1}(x)g(x)\right)\,(T-1)^{d'-1},$$}
 which implies that for $d-n>p$, i.e. $d>p+n$,
we have  $u_{d,p}(x)=0$, while for $d-n\le p< m$ or $d-n<p\le m$,
{\scriptsize
\begin{eqnarray*}u_{d,p}(x)&=& (-1)^{n(d-m)+d}\big(
(-1)^{n-q-1}{d'-1\choose m-p}F_{d'-1}(x)f(x)+
(-1)^{n-q}{d'-1\choose
m-p-1}G_{d'-1}(x)g(x)\big)\,\V(A)\,V(B)\\[2mm]
&=&(-1)^{n(d-m)+d+n-q-1}\big( {d'-1\choose m-p}F_{d'-1}(x)f(x)-
{d'-1\choose m-p-1}G_{d'-1}(x)g(x)\big)\,\V(A)\,V(B)
\end{eqnarray*}}
Therefore, by (\ref{th2.1}), for $n\le d\le m+n-1$ we have
{\scriptsize
\begin{eqnarray*}\operatorname{Sylv}^{p,q}(A,B;x)&=&(-1)^{q(m-p)+n(d-m)+d+n-q-1}
\big( {d'-1\choose m-p}F_{d'-1}(x)f(x)- {d'-1\choose
n-q}G_{d'-1}(x)g(x)\big),
\end{eqnarray*}}
since {\small $$ {d'-1\choose m-p-1}={d'-1\choose
d'-m+p}={d'-1\choose n-q}.
$$}
Finally for $d=m+n$, i.e. $p=m,q=n$ we have

$$u_d(x,T)=\V(A)\,\V(B)\,\operatorname{Sylv}^{m,n}(A,B;x)=
\V(A)\,\V(B)\,\Res(f,g)f(x)g(x)$$ which implies the claim.
The main theorem has been proved.
\end{pf}

%
%
%

\bigskip
\begin{ack}
T.~Krick would like to thank the Institute for Mathematics and its
Applications at Minneapolis and  the Department of Mathematics of
North Carolina State University where she was a guest during the
preparation of this note.
\end{ack}

\bigskip


\begin{thebibliography}{}
\bibitem[Ap\'ery and Jouanolou(2006)]{AJ06}
Ap\'ery, Fran\c cois ; Jouanolou, Jean-Pierre.
\newblock{\em R\'esultant et sous-r\'esultant :
le cas d'une variable avec exercices corrig\'es.\/}
Hermann, Paris 2006. 477 p.

\bibitem[{Borchardt(1860)}]{Bch60}
Borchardt, Carl Wilhelm.
\newblock{\em Uber eine Interpolationsformel f\"ur eine Art symmetrischer
Funktionen und \"uber deren Anwendung.\/}
\newblock Math. Abh. Akad. Wiss. zu Berlin (1860), 1--20.

\bibitem[{Borchardt(1878)}]{Bch78}
Borchardt, Carl Wilhelm.
\newblock{\em Zur Theorie der Elimination und
Kettenbruchentwicklung, Math.\/}
\newblock Abh. Akad. Wiss. zu Berlin (1878), 1--17.

\bibitem[{D'Andrea et al.(2007)}]{DHKS}
D'Andrea, Carlos; Hong, Hoon; Krick, Teresa; Szanto, Agnes.
\newblock{\em  An elementary proof of Sylvester's double sums for subresultants.\/}
\newblock  J. Symbolic Comput.  42  (2007),  no. 3, 290--297.

\bibitem[{Lascoux and Pragacz(2003)}]{LP}
Lascoux, Alain; Pragacz, Piotr. \newblock{\em Double Sylvester sums
for subresultants and multi-Schur functions.\/} \newblock J.
Symbolic Comput. 35 (2003),  no. 6, 689--710.

\bibitem[{Sylvester(1853)}]{sylv}
Sylvester, James Joseph.
\newblock {\em On a theory of syzygetic relations of two rational integral
functions, comprising anapplication to the theory of Sturm's
function and that of the greatest algebraical common measure.\/}
\newblock Philosophical Transactions of the Royal Society of London, Part III
(1853), 407--548.
\newblock Appears also in Collected Mathematical Papers
          of James Joseph Sylvester, Vol 1,
\newblock Chelsea Publishing Co. (1973), 429-586.

\bibitem[{Sylvester(1973)}]{sylv2}
Sylvester, James Joseph.
\newblock{\em On a generalization of the Lagrangian theorem of interpolation.\/}
\newblock Philosophical Magazine (1858).
\newblock Appears also in Collected Mathematical Papers
          of James Joseph Sylvester, Vol 1,
\newblock Chelsea Publishing Co. (1973), 645--646.
\end{thebibliography}
\end{document}